\documentclass[11pt]{amsart}
\usepackage{amsmath}
\usepackage{amsthm}
\usepackage{amssymb}
\newtheorem{thm}{Theorem}[section]

\newtheorem{lem}[thm]{Lemma}

\theoremstyle{definition}
\newtheorem{defn}[thm]{Definition}
\newtheorem{ex}[thm]{Example}
\newtheorem{rem}[thm]{Remark}
\newtheorem{con}[thm]{Construction}

\def\rk{\text{\rm rk}}
\def\deg{\text{\rm deg}}
\def\int{\text{\rm int}}
\def\sign{\text{\rm sign}}
\def\tf{{^2\!\phi}}
\def\ctf{{^2\!\phi}}
\def\Arf{\text{\rm Arf}}
\def\d{\delta}
\def\e{\varepsilon}
\def\E{E}
\def\Z{\mathbb Z}

\def\R{\mathbb R}
\def\C{\mathbb C}
\def\H{\mathbb H}
\def\P{\mathbb P}
\def\F{\mathbb F}
\def\f{\text{\rm F}}
\def\tF{{^2\F}}
\def\o{\mathfrak o}
\def\Tor{\text{\rm Tor}}
\def\Spin{\text{\rm Spin}}

\def\heading#1{\penalty-200\bigskip\noindent\sl#1\medskip\rm}

\begin{document}

\baselineskip=13.5 pt plus 1pt minus .5pt

\title{Canonical framings for $3$-manifolds}
\author{Rob Kirby and Paul Melvin}
\address{Department of Mathematics, University of California, Berkeley, CA
94720, USA}
\email{kirby@math.berkeley.edu}
\address{Department of Mathematics, Bryn Mawr College, Bryn Mawr, PA
19010, USA}
\email{pmelvin@brynmawr.edu}

\maketitle


\setcounter{section}{-1}
\section{Introduction}

A framing of an oriented trivial bundle is a homotopy class of sections of
the associated oriented frame bundle.  This paper is a study of the
framings of the tangent bundle $\tau_M$ of a smooth closed oriented
$3$-manifold $M$, often referred to simply as {\it framings} of
$M$.\footnote
{Recall that $\tau_M$ is always trivial; see \cite[p.148]{ms}, or
\cite[p.46]{kir} for an elementary geometric proof.}
We shall also discuss {\it stable framings} and {\it $2$-framings} of $M$,
that is framings of $\varepsilon^1\oplus\tau_M$ (where $\varepsilon^1$ is
an {\sl oriented} line bundle) and $2\tau_M = \tau_M \oplus \tau_M$.

The notion of a {\it canonical $2$-framing} of $M$ was introduced by
Atiyah \cite{atiyah}.  Motivated by Witten's paper \cite{witten}
generalizing the Jones polynomial to links in $M$, Atiyah observed that
Witten's invariant contained a phase factor specified by the choice of a
2-framing on $M$, and thus was an invariant of links in a 2-framed
3-manifold. Independent calculations by Reshetikhin and Turaev \cite{rt}
for a related invariant, defined from a framed link description of $M$,
did not appear to depend on a 2-framing.  As Atiyah noted, however, the
framed link description naturally gave a 2-framing of $M$, explained
further by Freed and Gompf in \cite[\S2]{fg}, and so  Reshetikhin and
Turaev were in fact calculating Witten's invariant for $M$ with this
framing times a phase factor depending on the difference between this
framing and the canonical 2-framing, i.e.\ Witten's invariant for $M$ with
its canonical 2-framing.

In this paper we give a leisurely exposition of framings, stable framings
and 2-framings of $M$, including some of the material in \cite{atiyah} and
\cite{fg}.  Our principal objective is to define the notion of a {\it
canonical (stable) framing} within each spin structure on $M$.  This is the
content of \S2.  The set of possible framings $\varphi$ for a given spin
structure form an affine space $\Z$, corresponding to $\pi_3(SO_3)$, and we
choose a canonical framing in this space by minimizing the absolute value
of the ``Hirzebruch defect" $h(\varphi)$ (defined in \S1).  More
generally, there are $\Z \oplus \Z = \pi_3(SO_4)$ possible stable framings
$\phi$, and here we must also minimize a certain ``degree" $d(\phi)$
associated with $\phi$.

The typical application may be to calculate the difference between a
naturally occurring framing and the canonical one, and this is carried out
in a number of instances.  In \S3 we consider framings on quotients of
$S^3$ by finite subgroups, where the calculations use signature defects
and the G-signature theorem, and on certain circle bundles over surfaces.
Natural framings also arise by restriction from framings on $4$-manifolds
bounded by $M$, and this situation is taken up in \S4.  In particular we
discuss surgery on an even framed link $L$ in $S^3$, which connects
with the work of Freed and Gompf.

Some of this work was done around the time of Atiyah's paper \cite{atiyah},
whose appearance removed our impetus to publish in a timely fashion.  We
wish to thank Selman Akbulut, Turgut Onder and our Turkish hosts for a
splendid conference in Gokova which provided the motivation to complete
this work.


\section{Preliminaries}

We work throughout in the smooth category.  It is assumed that the reader
is familiar with the elementary notions of handlebody theory, as treated in
Chapter I of \cite{kir}, and with the basic theory of fibre bundles and
characteristic classes, as presented for example in the beautiful books
of Steenrod \cite{st} and Milnor-Stasheff \cite{ms}.

Our orientation conventions are as follows.  It is always assumed
implicitly that the framings of an oriented manifold $M$ are consistent
with its orientation.  Reversal of orientations will be indicated with an
overbar.  In particular a framing or stable framing $\phi$ is reversed by
negating the last vector in each frame, producing $\bar\phi$.

We use the ``outward normal first" convention for compatibly orienting a
manifold and its boundary.  In particular, if $M$ is the oriented boundary
of an oriented $4$-manifold $W$, then the oriented bundles
$\varepsilon^1\oplus\tau_M$ and $\tau_W|M$ are naturally isomorphic by
identifying a framing $\nu$ of $\varepsilon^1$ with the outward pointing
normal.  This identification will be used implicitly when we discuss the
problem of extending stable framings of $M$ across $W$.

\heading{The degree and Hirzebruch defect}

Each framing $\varphi$ of a $3$-manifold $M$ can be identified with the
stable framing $\phi = \nu\oplus\varphi$, where $\nu$ is a framing of
$\varepsilon^1$ as above.  (Note the difference between the symbols
$\varphi$ for the ``honest" framing and $\phi$ for the stable framing.)  Of
course not all stable framings arise in this way.  The obstruction to
``destabilizing" a stable framing $\phi$ is measured by the {\it degree}
$d(\phi)$ of the map $M \to S^3$ which assigns to each point of $M$ the
position of the outward normal in the $3$-sphere determined $\phi$,
$$
d(\phi) = \deg(\,\nu\!:M\to S^3\,).
$$
Other properties of the degree are discussed in $\S2$.

In addition to the degree, the key invariant that will be brought to bear
on the study of framings is the {\it Hirzebruch defect}, defined for a
framing or stable framing $\phi$ of $M$ by
$$
h(\phi) = p_1(W,\phi) - 3\sigma(W)
$$
where $W$ is any compact oriented $4$-manifold bounded by $M$.  Here
$p_1(W,\phi)$ is the relative first Pontrjagin number of $W$ (explained
carefully in Appendix A), and $\sigma(W)$ is the signature of $W$.  It
follows from the Hirzebruch signature theorem and Novikov additivity of
the signature that this definition is independent of the choice of $W$.
For $2$-framings $\tf$ of $M$ we define
$
h(\tf) = p_1(W,\tf) - 6\sigma(W)
$
\cite{atiyah}.  These definitions are motivated and explained more
thoroughly in Appendix B, which includes a general discussion of signature
theorems and ``defect" invariants of $3$-manifolds.

It will be seen in the next section that the degree and Hirzebruch defect
form a complete set of integer invariants for the stable framings within
any given spin structure on $M$, and serve to identify this set of framings
with an affine lattice in the the $dh$-plane.  The {\sl canonical}
framing(s) will then be defined as the one(s) closest to the origin in
this plane, with respect to a suitable norm.  Before embarking on this
program, however, we identify two well known elements of $\pi_3(SO_4)$ that
will be used frequently in what follows.

\heading{The generators $\rho$ and $\sigma$}

The study of framings leads naturally to obstruction theory with
coefficients in the homotopy groups $\pi_3(SO_n)$ and $\pi_3(SU_n)$.
Recall that $\pi_3(SO_4) \cong \Z\oplus\Z$, $\pi_3(SO_n) \cong
\Z$ for all other $n\ge3$, and $\pi_3(SU_n) \cong \Z$ for all $n\ge2$.
Following Steenrod \cite{st}, we give explicit generators for all of these
groups.

View $S^3$ as the unit sphere in the quaternions $\H$ (oriented by the
ordered basis $1,i,j,k$) and $SO_4$ as the rotation group of $\H$.  Then
the maps $\rho$ and $\sigma\!:S^3\to SO_4$ defined by
$$
\rho(q)x = qxq^{-1} \qquad \sigma(q)x = qx
$$
represent generators of $\pi_3(SO_4)$ \cite[p.117]{st}.  (By abuse of
notation, we shall not distinguish between these maps and their homotopy
classes.)  This identifies $\pi_3(SO_4)$ with $\Z\oplus\Z$, where $\rho$,
$\sigma$ correspond to the standard ordered basis.

By restricting to the subspace $\P$ of pure (imaginary) quaternions,
$\rho$ also represents an element (in fact a generator) of $\pi_3(SO_3)$,
and these two $\rho$'s correspond under the natural map $\pi_3(SO_3)\to
\pi_3(SO_4)$ induced by the inclusion $SO_3\subset SO_4$.  Furthermore,
$\sigma \in
\pi_3(SO_4)$ is carried under the natural maps $\pi_3(SO_4) \to
\pi_3(SO_5) \to \cdots \to \pi_3(SO)$ onto generators, also denoted by
$\sigma$, of each of the subsequent groups; the first of these maps also
carries $\rho$ to $2\sigma$.  These facts are all proved in Steenrod, and
serve to identify all of these groups with $\Z$.

Similarly $\sigma$ represents a generator of $\pi_3(SU_2)$ by identifying
$\H$ with $\C^2$.  (Note that for $\sigma$ to induce a {\sl unitary}
action, the vector $(u,v)$ in $\C^2$ must be identified with the quaternion
$u+jv$ rather than $u+vj$; see also Remark \ref{rem:p} below.)
This generator is then carried onto generators $\sigma$ for all the
subsequent groups in the natural sequence $\pi_3(SU_2) \to \pi_3(SU_3) \to
\cdots \to \pi_3(SU)$, thereby identifying these groups with $\Z$.

Finally observe that the tautological map
$$
\iota\!:\pi_3(SO) \to \pi_3(SU)
$$
is multiplication by $-2$, that is $\iota(\sigma) = -2\sigma$.  Indeed,
this map was classically computed up to sign (see e.g.\ the proof of Lemma
2 in \cite{mke}), and an explicit formula giving the sign was written down
by John Hughes in his thesis \cite[\S24]{hughes}, essentially as follows:

Fix a unit quaternion $q = a + bi + cj + dk = u+jv$, with $u = a+bi$ and $v
= c-di$, and let $P$ and $Q$ denote the matrices in $SO_4$ and $SU_2$
corresponding to $\sigma(q)$,
$$
P = \pmatrix
a & -b & -c & -d \\
b & a & -d & c \\
c & d & a & -b \\
d & -c & b & a
\endpmatrix
\qquad\text{and}\qquad
Q = \pmatrix
u & -\bar v \\
v & \bar u
\endpmatrix.
$$
Now observe that $P$ is conjugate in $SU_4$ to the block sum $T = Q^{-1}
\oplus Q^{-1}$ of two copies of the inverse of $Q$.  Explicitly $RPR^{-1}
= T$ for
$$
R = {1\over\sqrt2}
\pmatrix
0 & 0 & 1 & -i \\
1 & i & 0 & 0 \\
i & 1 & 0 & 0 \\
0 & 0 & -i & 1
\endpmatrix.
$$
Since $SU_4$ is connected, it follows that the maps $\pi$ and $\tau \!:
S^3\to SU_4$ given by $\pi(q) = P$ and $\tau(q) = T$ are homotopic, and
evidently $\pi$ represents $\iota(\sigma)$ and $\tau$ represents
$-2\sigma$.

\heading{The geometry of $\rho$ and $\sigma$}

There is a simple geometric description for the maps $\rho$ and $\sigma
\in \pi_3(SO_4)$ if one views $S^3$ as the unit $3$-ball $B^3$ in $\P$ with
its boundary collapsed to a point.  (The great half circles in $S^3$ from
$-1$ to $1$ correspond to the radii of $B^3$.)   In particular $\rho(0) =
I$ and $\sigma(0) = -I$.  If $x\in B^3$ is nonzero, then $\rho(x)$ is the
rotation of $\P$ about $x$ by $2\pi|x|$ radians, while $\sigma(x)$ is the
simultaneous rotation by $\pi|x|$ radians of the oriented plane $P_x$
(normal to $x$ in $\P$) and its complement $Q_x$ (with ordered basis $1$,
$x$).  In other words as the diameter through $x$ is traversed, $\rho$
rotates $P_x$ two full {\sl right}-handed turns while fixing $Q_x$, and
$\sigma$ rotates both planes by one full right-handed turn.


\section{Canonical framings}

Let $M$ be a closed connected oriented $3$-manifold.  We wish to study the
framings $\varphi$, stable framings $\phi$, and $2$-framings $\tf$ of
$M$.  To organize this study, it is convenient to fix a spin structure
$\Sigma$ on $M$, which can be viewed as a framing of $\tau_M$ over the
$2$-skeleton of $M$ (see \cite{mil} and the discussion in Chapter
IV of \cite{kir}).

\heading{Framings extending $\Sigma$}

Consider the set $\f_\Sigma$ of framings of $M$ which are compatible with
$\Sigma$.  The difference of two such framings is specified by an element
of $H^3(M;\pi_3(SO_3)) = \Z$, by obstruction theory.  In other words
$\f_\Sigma$ is an affine space with translation group $\pi_3(SO_3) = \Z$.
The action of this group can be visualized geometrically using the
description in \S1 of its generator $\rho$.  In particular $\rho$ acts on a
framing $\varphi$, producing a new framing $\varphi + \rho$, by rotating
the frames along each diameter of a small $3$-ball in $M$ by two full
twists.  It will be seen below (Lemma \ref{lem:p}a) that this corresponds
to a shift by $4$ in the affine space $\f_\Sigma$.

\heading{Canonical framings}

As an affine space, $\f_\Sigma$ has no apriori choice of basepoint.  Our
goal is to pick such a basepoint, i.e.\  a preferred or {\sl canonical}
framing in the given spin structure on $M$.  We shall not attempt to pick a
canonical framing within the set $\f$ of {\sl all} framings on $M$, or what
amounts to the same thing, a canonical spin structure on $M$.  Indeed we
do not know how to make such a choice in general.  Recall that by
obstruction theory, as above, the spin structures on $M$ form an affine
space over $\Z_2$ with translation group $H^1(M;\Z_2)$, and so $\f =
\cup\f_\Sigma$ is (non-canonically) isomorphic to $H^1(M;\Z_2) \oplus \Z$.

\begin{defn}\label{defn:cf}
A framing $\varphi$ of the $3$-manifold $M$ is {\it canonical} for the spin
structure $\Sigma$ if it is compatible with $\Sigma$ and $|h(\varphi)|
\le |h(\psi)|$ for all other framings $\psi$ which are compatible
with $\Sigma$.  In other words, $\varphi$ is a minimum for the invariant
$|h|$ on $F_\Sigma$.
\end{defn}

Since $\rho$ translates $\f_\Sigma$ by $4$ (Lemma \ref{lem:p}a below),
this minimum is at most $2$.  If it is $0$ or $1$, then the canonical
framing is unique (in $\f_\Sigma$) and in fact minimizes $|h|$ globally (in
$\f$).  This is the case for ``most" spin structures on most manifolds.
The only exceptions are the spin structures for which $\mu\equiv2\rk \,
H_1(M;\Z_2) \pmod4$ (e.g.\ the unique spin structure on a homology
sphere).  In these cases there are two canonical framings with $h=\pm2$;
one could of course select the positive one, but we have chosen not to out
of respect for our left-handed colleagues.

Note that if one were to define a canonical framing to be one that
minimized $|h|$ globally, then some spin structures might not have {\sl
any} canonical framings; indeed examples are easily given of manifolds in
which the minimum of $|h|$ varies with the spin structure (e.g.\ the
connected sum of two copies of real projective $3$-space).

To justify the preceding remarks, it is necessary to analyze the behavior
of the defect $h(\varphi)$ under the action of $\pi_3(SO_3)$ on the set of
framings.  This is not difficult, and is discussed below in the more
general context of stable framings.

\heading{Stable framings extending $\Sigma$}

Let $\F$ be the set of all stable framings of $M$, and $\F_\Sigma$
be the subset of those which extend the spin structure
$\nu\oplus\Sigma$ on $\varepsilon^1\oplus\tau_M$, where $\nu$ is a framing
of $\varepsilon^1$ (restricted to the $2$-skeleton of $M$).  As above,
$\F_\Sigma$ is an affine space with translation group $\pi_3(SO_4) =
\Z\oplus\Z$, and $\F \cong H^1(M;\Z_2) \oplus \Z \oplus
\Z$ (non-canonically). The action of $\pi_3(SO_4)$ on $\F$ can be
understood as before using the geometric descriptions of $\rho$ and
$\sigma$ (with $\nu$ playing the role of $1$):  $\phi + m\rho + n\sigma$
is obtained from $\phi$ by rotating the frames along each oriented
diameter of a $3$-ball by $2m+n$ full twists in the normal plane,
perpendicular to the diameter, and by $n$ full twists in the ``conormal"
plane, spanned by $\nu$ and the diameter.

\heading{The degree}

Observe that the ``honest" framings $\f_\Sigma$ can be viewed as a subset
of the stable framings $\F_\Sigma$ by identifying $\varphi$ with $\phi =
\nu\oplus\varphi$.  Thus we will freely use $\varphi$ and $\phi$
interchangeably with the understanding that this identification is to be
made.  There is a simple invariant that detects whether a stable framing
$\phi$ corresponds to an honest framing in this way, namely the {\it
degree} of
$\phi$
$$
d(\phi) = \deg(\,\nu\!:M\to S^3\,).
$$
Here $S^3$ is the unit sphere in $\phi$, and the map $\nu\!:M\to S^3$ is
defined in the obvious way: $\nu$ is the framing of $\varepsilon^1$, which
at each point of $M$ determines a point in $S^3$.

The degree satisfies the following properties.

\begin{thm}\label{thm:d} Let $\phi$ be a stable framing of $M$.  Then
$d(\phi) = 0$ if and only if $\phi$ is of the form $\nu\oplus\varphi$ for
some honest framing $\varphi$.  Furthermore,
\begin{enumerate}
\item[a)] {\rm(action)} \ $d(\phi+\rho) = d(\phi)$ and
$d(\phi+\sigma) = d(\phi)-1$.
\item[b)] {\rm(boundary)} \ If $\phi$ extends to a framing of a
compact $4$-manifold $W$ bounded by $M$, then $d(\phi) = \chi(W)$, where
$\chi$ is the Euler characteristic.
\item[c)] {\rm(covering)} \ If $(\tilde M,\tilde\phi) \to (M,\phi)$ is an
$r$-fold covering map with compatible stable framings, then $d(\tilde\phi)
= r\,d(\phi)$.
\item[d)] {\rm(orientation)} \ $d(\bar\phi) = d(\phi)$, where \
$\bar{ }$ \ denotes orientation reversal.
\end{enumerate}
\end{thm}

\begin{proof}
The first statement is immediate from the fact that the degree of a map
$M\to S^3$ is zero if and only if the map is homotopic to a constant
(Hopf's Theorem).  Property a) can be verified using the geometric
description of the action of $\rho$ and $\sigma$ on $\F$, and b) follows
from the fact that the Euler characteristic is equal to the obstruction to
extending $\nu$ across $W$, since the latter can be identified with
$d(\phi)$ by the homological invariance of degree.  (Thus $d(\phi)$ can be
interpreted as the relative Euler class of $(W,\phi)$.)  Property c) is
obvious.  For d) we have $d(\bar\phi) = \deg(\,\nu\!:\bar M\to\bar S^3\,) =
d(\phi)$.
\end{proof}

\heading{The Pontrjagin number}

Next we investigate the relative first Pontrjagin number $p_1(W,\phi)$,
where $W$ is a compact $4$-manifold bounded by $M$.  This is a key
ingredient in the definition of the Hirzebruch defect $h(\phi)$, and was
defined in Appendix A as the obstruction to extending $\phi$ (with the last
vector dropped) across the complexified tangent bundle of $W$.

If $W$ has a spin structure compatible with $\phi$, then
$p_1$ can be defined in a simpler way without complexifying $\tau_W$.  For
in this case $\phi$ extends to a framing of $W$ in the complement of a
point.  Indeed the spin structure gives an extension over the relative
$2$-skeleton of $(W,\partial W)$, and this further extends over the
$3$-skeleton since  $\pi_2(SO_4)=0$.  Now the obstruction to extending
$\phi$ across the last point is an element $\theta$ in $\pi_3(SO_4) =
\Z\oplus\Z$.  This element defines a 4-plane bundle $\xi_{\theta}$ over
$S^4$ whose first Pontrjagin number $p_1(\xi_\theta)$ evidently coincides
with $p_1(W,\phi)$.  Alternatively, one may consider the image $\o$ of
$\theta$ in $\pi_3(SO) = \Z$, the ``stable obstruction" to extending
$\phi$, and it is easily verified that
$$
p_1(W,\phi) = 2\o
$$
since $\pi_3(SO) \to \pi_3(SU)$ is multiplication by $-2$ and $p_1=-c_2$.

The relative first Pontrjagin number satisfies the following properties.

\penalty-800
\begin{lem}\label{lem:p}
 Let $\phi$ be a stable framing of $M = \partial W$.  Then
\begin{enumerate}
\item[a)] {\rm(action)} \ $p_1(W,\phi+\rho) = p_1(W,\phi)+4$ and
$p_1(W,\phi+\sigma) = p_1(W,\phi)+2$.
\item[b)] {\rm(boundary)} \ If $\phi$ is compatible with a spin structure
on $W$, then $\phi$ extends to a framing of $W$ if and only if
$p_1(W,\phi) = 0$ and $d(\phi) = \chi(W)$.
\item[c)] {\rm(covering)} \ If $(\tilde W,\tilde\phi) \to (W,\phi)$
is an $r$-fold covering map with compatible stable framings on the
boundary, then $p_1(\tilde W,\tilde\phi) = rp_1(W,\phi)$.
\item[d)] {\rm(orientation)} \ $p_1(\bar W,\bar\phi) = -p_1(W,\phi)$, where
\
$\bar{ }$ \ denotes orientation reversal.
\end{enumerate}
\end{lem}

\begin{proof}
The action of any $\theta \in \pi_3(SO_4)$ on $\phi$ is local, only
changing $\phi$ in a $3$-ball in $M$, and so  $p_1(W,\phi+\theta)
= p_1(W,\phi) + p_1(B^4, \phi_0 + \theta)$, where $\phi_0$ is the
restriction to $S^3$ of the unique framing of $B^4$.  But $p_1(B^4, \phi_0
+ \theta) = 2\o$, where $\o$ is the image of $\theta$ in $\pi_3(SO) = \Z$.
Evidently $\o = 1$ or $2$ for $\theta = \sigma$ or $\rho$, repectively, and
the property a) follows.

The forward implication in b) is immediate from Theorem \ref{thm:d}b
and the definition of $p_1$.  For the reverse implication, suppose that
$\phi\in\F_\Sigma$.  Observe that any spin structure on $W$ which extends
$\Sigma$ can be further extended to a framing of $W$, by obstruction
theory.  Let $\psi$ denote the restriction of any such framing to $M$.
Then $\psi = \phi + m\rho + n\sigma$, for some $m$ and $n$, since both
stable framings extend $\Sigma$.  By hypothesis $p_1(W,\phi) = p_1(W,\psi)
= 0$, and so $2m+n = 0$ by a).  Also $d(\phi) = d(\psi) = \chi(W)$ and
so $n = 0$ by Theorem \ref{thm:d}a.  Thus $m=n=0$, and so
$\phi=\psi$.  This establishes the b).  Properties c) and d) are immediate
from the definition of $p_1$ as an obstruction (see Appendix A).
\end{proof}

\begin{rem}\label{rem:p}
Lemma \ref{lem:p}a can also be established using the formula $p_1(B^4,
\phi_0+\theta) = p_1(\xi_\theta)$ (see the discussion preceding the
lemma).  Indeed the {\sl oriented} bundle $\bar\xi_\sigma$ can be given
the structure of a complex $2$-plane bundle $\omega$ over $S^4$, by right
multiplication, and so
$$
p_1(\xi_\sigma ) = c_1^2(\omega)-2c_2(\omega) = -2e(\bar\xi_\sigma) = 2.
$$
Note that the natural orientation on $\xi_\sigma$ coming from $\H$ is
inconsistent with the orientation arising from the complex structure.
This is the reason for using the oppositely oriented bundle
$\bar\xi_\sigma$, with Euler class $-1$, to compute $c_2(\omega)$ (cf.\
the computations in \cite[p.673]{dw} and \cite[p.43]{kir} where this sign
is overlooked). Now since $\tau_{S^4} = \xi_{2\sigma - \rho}$
(\cite[\S23.6,27.3]{st}) and $p_1(\tau_{S^4}) =
p_1(\varepsilon^1\oplus\tau_{S^4}) = p_1(\varepsilon^5) = 0$, it follows
that $p_1(\xi_{\rho}) = 4$.
\end{rem}

\heading{The Hirzebruch defect}

Using the previous lemma, together with some facts about the signature and
signature defects (see Appendix B), we deduce the following properties of
the Hirzebruch defect $h(\phi)$.  (Recall from \S1 that $h(\phi) =
p_1(W,\phi) - 3\sigma(W)$ for any compact oriented $4$-manifold bounded
by $M$.)

\penalty-800
\begin{thm}\label{thm:h}
Let $\phi$ be a stable framing of $M$.  Then
\begin{enumerate}
\item[a)] {\rm(action)} \ $h(\phi+\rho) = h(\phi)+4$ and
$h(\phi+\sigma) = h(\phi) + 2$.
\item[b)] {\rm(boundary)} \ If $\phi$ extends to a framing of a compact
$4$-manifold $W$ bounded by $M$, then $h(\phi) = -3\sigma(W)$, where
$\sigma(W)$ is the signature of $W$.
\item[c)] {\rm(covering)} \ If $\pi\!:(\tilde M,\tilde\phi) \to (M,\phi)$
is an $r$-fold cover with compatible stable framings, then
$h(\phi) = (h(\tilde\phi) - 3\sigma(\pi))/r$, where
$\sigma(\pi)$ is the signature defect of $\pi$.
\item[d)] {\rm(orientation)} \ $h(\bar\phi) = -h(\phi)$, where \
$\bar{ }$ \ denotes orientation reversal.
\end{enumerate}
\end{thm}

\begin{proof}
Properties a) and b) are immediate from the analogous properties for $p_1$
in Lemma 2.3.  Property c) is proved in Appendix B, Lemma 1.  The last
property follows from Lemma \ref{lem:p}d and the fact that $\sigma(\bar W)
= -\sigma(W)$.
\end{proof}

\heading{The total defect}

The preceding results give a complete picture of the affine space
$\F_\Sigma$ of stable framings which extend the spin structure $\Sigma$ on
$M$.  Indeed Theorems \ref{thm:d}a and \ref{thm:h}a show that $d$ and
$h$ together give an embedding
$
H:\F_\Sigma \hookrightarrow \Z\oplus\Z
$
into the $dh$-plane, defined by
$$
H(\phi) = (d(\phi),h(\phi)).
$$
This embedding will be called the {\it total defect}.  The (vertical)
$h$-axis corresponds to the honest framings $\varphi$, and in line with
the natural inclusion $\f_\Sigma \subset \F_\Sigma$ we often write
$H(\varphi)$ for $H(\nu\oplus\varphi)$.  Orientation reversal
$\phi \mapsto \bar\phi$ corresponds to ``conjugation" $(d,h) \mapsto
(d,-h)$ (by \ref{thm:d}d and \ref{thm:h}d).

>From the action of $\pi_3(SO_4)$ on the $dh$-plane (shown in Figure 1a),
the image of $\F_\Sigma$ is seen to be an affine lattice of index $4$ in
$\Z\oplus \Z$.  In particular it is a coset of the subgroup $\Lambda_0$
generated by $(0,4)$ and $(-1,2)$ (shown in Figure 1b), namely one of the
four affine lattices $\Lambda_k = \Lambda_0 + (0,k)$ for $k \in \Z_4$.

\begin{center}
\begin{picture}(250,120)

\put(0,0)
{\begin{picture}(100,160)
\multiput(10,20)(0,20){5}{\thinlines\line(1,0){80}}
\multiput(20,10)(20,0){4}{\thinlines\line(0,1){100}}
\put(45,30){$\sigma$}
\put(64,67){$\rho$}
\put(60,20){\thicklines\vector(0,1){80}}
\put(60,20){\thicklines\vector(-1,2){20}}
\put(-5,-5){a) the action of $\pi_3(SO_4)$}
\end{picture}}

\put(150,0)
{\begin{picture}(100,160)
\multiput(10,20)(0,20){2}{\thinlines\line(1,0){80}}
\put(10,60){\thinlines\vector(1,0){80}}
\put(95,58){$d$}
\multiput(10,80)(0,20){2}{\thinlines\line(1,0){80}}
\put(20,10){\thinlines\line(0,1){100}}
\put(40,10){\thinlines\vector(0,1){100}}
\put(36,114){$h$}
\multiput(60,10)(20,0){2}{\thinlines\line(0,1){100}}
\multiput(20,20)(40,0){2}{\circle*{4}}
\multiput(40,60)(40,0){2}{\circle*{4}}
\multiput(20,100)(40,0){2}{\circle*{4}}
\put(15,-5){b) the lattice $\Lambda_0$}
\end{picture}}

\end{picture}
\end{center}

\medskip
\begin{center}
Figure 1: the $dh$-plane
\end{center}
\medskip

\noindent
To determine which one, consider the $\Z_4$-valued invariant $\lambda =
2d+h \pmod4$ of the spin structure $\Sigma$ on $M$, that is
$$
\lambda(\Sigma) \ = \ 2d(\phi) + h(\phi) \!\!\!\pmod 4
$$
for any stable framing $\phi$ in $\F_\Sigma$.  Then
$
H(\F_\Sigma) = \Lambda_{\lambda(\Sigma)},
$
and the next result expresses the invariant $\lambda(\Sigma)$, and its mod
$2$ reduction, in terms of the mu invariant $\mu(\Sigma)$ and homological
invariants of $M$.  (Recall that $\mu(\Sigma) = \sigma(W) \pmod{16}$ for
any compact $4$-manifold $W$ bounded by $M$ over which $\Sigma$ extends;
this is well defined by Rohlin's theorem on the signature of closed spin
$4$-manifolds \cite{roh} \cite[pp.64--65]{kir}.)

\begin{thm}\label{thm:l}  Let $\Sigma$ be a spin structure on $M$.  Then
\begin{align}
\lambda(\Sigma) \ &\equiv \ 2(1+r(M))+\mu(\Sigma) \!\!\!\pmod4 \notag\\
&\equiv \ s(M) \!\!\!\pmod2 \notag
\end{align}
where $r(M) = \rk(H_1(M) \otimes \Z_2)$ and $s(M) = \rk(\Tor H_1(M)
\otimes \Z_2)$.
\end{thm}

\noindent Note that $s(M)$ is the number of $2$-primary summands in
$H_1(M)$, and so $r(M)-s(M)$ is the first Betti number $b_1(M) =
\rk(H_1(M)).$

\begin{proof}
Choose any simply-connected spin $4$-manifold $W$ with spin boundary
$(M,\Sigma)$, for example constructed by attaching $2$-handles to $B^4$
along an even framed link \cite{kap}.  To compute $\lambda(\Sigma)$, we use
the restriction $\phi$ to $M$ of the unique framing of $W$.

First observe that $$\chi(W) \equiv 1+r(M) \pmod2.$$  Indeed $\chi(W) = 1 +
b_2(W)$, where $b_2(W) = \rk(H_2(W))$ can be expressed as the sum of
the nullity and rank of the mod $2$ intersection form on $W$.  But the
nullity of this form is equal to $r(M)$ since the intersection matrix is a
presentation matrix for $H_1(M;\Z_2)$ (by Poincar\'e duality), and the
rank is even since $W$ is spin (whence the nonsingular part of the form is
a sum of two-dimensional hyperbolic forms).

Now using Theorems \ref{thm:d}b and \ref{thm:h}b, we see that
$\lambda(\Sigma)$ is equal to
$$
2d(\phi) + h(\phi) \ \equiv \ 2\chi(W)-3\sigma(W) \ \equiv \ 2(1+r(M))
+ \mu(\Sigma) \pmod4.
$$
Reducing mod $2$ gives $\mu(\Sigma) \equiv \sigma(W) \equiv b_2(W) - b_1(M)
\equiv r(M) - b_1(M) \equiv s(M)$.
\end{proof}

\heading{Canonical stable framings}

The description of $\F_\Sigma$ above suggests the following generalization
of the notion of canonical framings (Definition \ref{defn:cf}) using the
norm on the $dh$-plane given by  $|(d,h)| = 2|d|+|h|$.

\begin{defn}\label{defn:csf}
A stable framing $\phi$ of the $3$-manifold $M$ is {\it canonical} for
the spin structure $\Sigma$ if it is compatible with $\Sigma$ and
$|H(\phi)| \le |H(\psi)|$ for all other stable framings $\psi$ which are
compatible with $\Sigma$.  In other words, $\phi$ is a minimum for the
invariant $2|d|+|h|$ on $\F_\Sigma$.
\end{defn}

It follows from this definition that $\F_\Sigma$ has a {\sl unique} stable
canonical framing when $\lambda(\Sigma) \not\equiv 2$, corresponding to
the point $(0,0)$ or $(0,\pm1)$ in the $dh$-plane, according to
whether $\lambda(\Sigma) \equiv 0$ or $\pm1$ (see Figure 2a-c).  These are
all honest framings, and so give the unique canonical framings in
$F_\Sigma$ as well.  If $\lambda(\Sigma) \equiv 2$ (e.g.\ for the unique
spin structure on any homology sphere) then there are {\sl four} canonical
framings, corresponding to the points $(\pm1,0)$ and $(0,\pm2)$ (see
Figure 2d).  The last two are honest framings.

\begin{center}
\begin{picture}(365,120)

\put(0,0)
{\begin{picture}(80,160)
\multiput(10,20)(0,20){2}{\thinlines\line(1,0){60}}
\put(10,60){\thinlines\vector(1,0){60}}
\put(75,58){$d$}
\multiput(10,80)(0,20){2}{\thinlines\line(1,0){60}}
\put(20,10){\thinlines\line(0,1){100}}
\put(40,10){\thinlines\vector(0,1){100}}
\put(36,114){$h$}
\put(60,10){\thinlines\line(0,1){100}}
\put(40,60){\circle*{4}}
\put(14,-5){a) $\lambda(\Sigma) = 0$}
\end{picture}}

\put(95,0)
{\begin{picture}(80,160)
\multiput(10,20)(0,20){2}{\thinlines\line(1,0){60}}
\put(10,60){\thinlines\vector(1,0){60}}
\put(75,58){$d$}
\multiput(10,80)(0,20){2}{\thinlines\line(1,0){60}}
\put(20,10){\thinlines\line(0,1){100}}
\put(40,10){\thinlines\vector(0,1){100}}
\put(36,114){$h$}
\put(60,10){\thinlines\line(0,1){100}}
\put(40,80){\circle*{4}}
\put(14,-5){b) $\lambda(\Sigma) = 1$}
\end{picture}}

\put(190,0)
{\begin{picture}(80,160)
\multiput(10,20)(0,20){2}{\thinlines\line(1,0){60}}
\put(10,60){\thinlines\vector(1,0){60}}
\put(75,58){$d$}
\multiput(10,80)(0,20){2}{\thinlines\line(1,0){60}}
\put(20,10){\thinlines\line(0,1){100}}
\put(40,10){\thinlines\vector(0,1){100}}
\put(36,114){$h$}
\put(60,10){\thinlines\line(0,1){100}}
\put(40,40){\circle*{4}}
\put(12,-5){c) $\lambda(\Sigma) = -1$}
\end{picture}}

\put(285,0)
{\begin{picture}(80,160)
\multiput(10,20)(0,20){2}{\thinlines\line(1,0){60}}
\put(10,60){\thinlines\vector(1,0){60}}
\put(75,58){$d$}
\multiput(10,80)(0,20){2}{\thinlines\line(1,0){60}}
\put(20,10){\thinlines\line(0,1){100}}
\put(40,10){\thinlines\vector(0,1){100}}
\put(36,114){$h$}
\put(60,10){\thinlines\line(0,1){100}}
\put(40,100){\circle*{4}}
\put(20,60){\circle*{4}}
\put(60,60){\circle*{4}}
\put(40,20){\circle*{4}}
\put(12,-5){d) $\lambda(\Sigma) = \pm2$}
\end{picture}}

\end{picture}
\end{center}

\medskip
\begin{center}
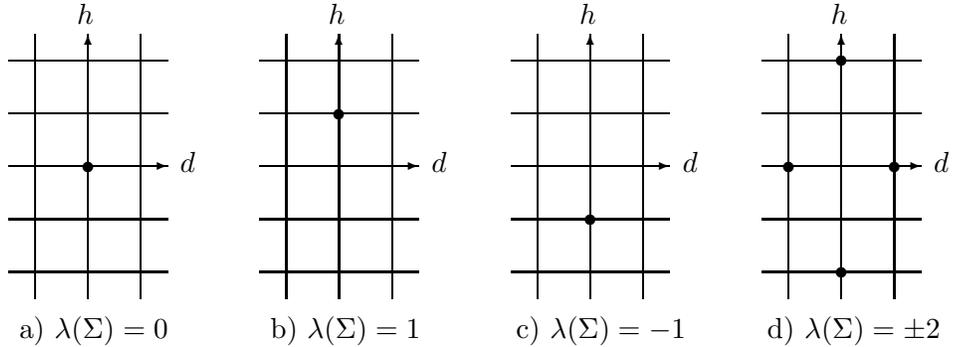
Figure 2: canonical stable framing(s) extending $\Sigma$
\end{center}
\medskip

\begin{rem}
Canonical framings can be constructed from any given stable
framing $\phi$ by adding suitable multiples of $\rho$ and $\sigma$.  For
example, if $\phi$ has total defect $(d,h)$, then the canonical framing
with total defect $(0,\lambda)$, where $\lambda \equiv 2d+h\pmod4$ and
$|\lambda|\le2$, is given by
$$
\phi + d\,\sigma - \textstyle{1\over4}(2d+h-\lambda)\,\rho.
$$
If $\phi$ can be described geometrically, then so can the canonical
framing, using the local geometric picture for the action of $\rho$ and
$\sigma$ discussed above.
\end{rem}

\heading{Examples of canonical framings}

\begin{ex}\label{ex:lie} (Lie groups) \ Consider the three $3$-dimensional
compact connected Lie groups -- the $3$-sphere $S^3$, the rotation group
$SO_3$, and the $3$-torus $T^3$ -- which have natural framings arising
from their group structure.  In particular $S^3$ and $SO_3$ each have two
Lie framings
$\varphi_\pm$, obtained by left or right multiplication from a fixed frame
at the identity, where the $+$ sign corresponds to left
multiplication.\footnote
{The reason for this sign convention in $S^3$ is
that {\sl left} multiplication yields the {\sl right}-handed Hopf framing,
in which the integral curves of any vectorfield in the framing are Hopf
circles with $+1$ pairwise linking, along which the other two vectorfields
spin once in a positive sense.  This is easily seen from the geometric
description for $\sigma$ in \S1.  The framing $\varphi_+$ of
$SO_3$ is likewise intrinsically right-handed, since it can be viewed as
the quotient of the corresponding framing of $S^3$ under the projection
$S^3 \to SO_3 = C_2\backslash S^3$.  Similar statements apply to the {\sl
left}-handed Hopf framing $\varphi_-$ on $S^3$ and its quotient on
$SO_3$.}   Since $T^3$ is abelian, it has only one Lie framing
$\varphi_1$.

Now $S^3$ has a unique spin structure, while $SO_3$ and $T^3$ have
$|H^1(SO_3;\Z_2)| = 2$ and $|H^1(T^3;\Z_2)| = 8$ spin structures,
respectively.  We will show how to construct {\sl all} the canonical
framings in these spin structures, and in particular show that the Lie
framings are canonical in theirs.  To emphasize the underlying manifold
$M$ of a framing $\varphi$, we will sometimes write $H(M,\varphi)$ for the
total defect $H(\varphi)$.

\smallskip
{\bf a)} For $S^3$, the Lie framings $\varphi_\pm$ have
Hirzebruch defect $\pm2$.  This can be seen in a variety of ways.  For
example, it is obvious that the obstruction to extending $\varphi_+$ over
$B^4$ is $\sigma$, and so $p_1(B^4,\varphi_+) = p_1(\xi_\sigma) = 2$ (as
shown in Remark \ref{rem:p}).  Thus $h(\varphi_+) = 2$
since $\sigma(B^4) = 0$.  Now observe that $\varphi_- = \varphi_+ - \rho$,
since right multiplication by $q$ is the composition of left
multiplication by $q$ with conjugation by $\bar q$ (i.e.\ $\bar q(qx)q =
xq$) and conjugation by $\bar q$ represents $-\rho \in \pi_3(SO_4)$.  It
follows that $h(\varphi_-) = -2$, by \ref{thm:h}a.

The same result can be obtained another way, using the ``canonical" stable
framing $\d$ of $S^3$ which is the restriction of the unique
framing of $B^4$.  Clearly $\varphi_+ = \d + \sigma.$  (This
provides a sign check for Theorem \ref{thm:d}a, that adding $\sigma$
lowers the degree by one: $\d$ has degree $1$ since $\chi(B^4) =
1$, while the honest framing $\varphi_+$ has degree $0$.)  Thus
$h(\varphi_+) = h(\d+\sigma) = 2$ by Theorem \ref{thm:h}a, since
$h(\d) = 0$ by \ref{thm:h}b.  Now $h(\varphi_-) = -2$ can be
deduced as above, or using Theorem \ref{thm:h}d and the observation that
$(S^3,\varphi_-) = (\bar S^3, \bar\varphi_+)$ (indeed any orientation
reversing automorphism of $S^3$ induces an orientation preserving
diffeomorphism $S^3 \to \bar S^3$ which identifies $\varphi_-$ with
$\bar\varphi_+$, up to homotopy).

Summarizing, we have computed the total defects $H = (d,h)$ of the Lie
framings $\varphi_\pm$ of $S^3$, and along the way the canonical stable
framing $\d$ coming from $B^4$, to be
$$
H(S^3,\varphi_\pm) = (0,\pm2) \qquad\text{and}\qquad H(S^3,\d) =
(1,0).
$$
Therefore these (stable) framings represent, by definition, three out of
four of the canonical stable framings in the unique spin structure on
$S^3$.  The fourth, represented by the point $(-1,0)$ in the $dh$-plane,
can also be constructed in a natural way as the restriction $\d_-$ of
{\sl any} framing of $B^4 \# S^1\times S^3$.  Indeed this $4$-manifold has
Euler characteristic $-1$ and signature $0$, and so $H(S^3,\d_-) =
(-1,0)$ (by Theorems \ref{thm:d}b and \ref{thm:h}b).

\smallskip
{\bf b)} For $SO_3$, note that $(S^3,\varphi_\pm)$ double covers
$(SO_3, \varphi_\pm)$ with zero signature defect (See Example 3 in Appendix
B).  It follows from Theorem \ref{thm:h}c that the Lie framings
$\varphi_\pm$ on $SO_3$ have Hirzebruch defect $\pm1$, and so total defects
$$
H(SO_3,\varphi_\pm) = (0,\pm1).
$$
Noting that the two spin structures $\Sigma_\pm$ on $SO_3$ have
$\mu$-invariants $\pm1$, it follows from Theorem \ref{thm:l} that
$\varphi_\pm \in \F_{\Sigma_\pm}$, and so we have identified the unique
canonical framings in both spin structures.  (These spin structures are
equivalent under any orientation reversing automorphism of $SO_3$, and in
fact $(SO_3,\varphi_-) = (\bar{SO}_3,\bar\varphi_+).$)

\smallskip
{\bf c)} For $T^3$, observe that the Lie framing $\varphi_1$ is
amphicheiral, i.e.\ $(T^3,\varphi_1) = (\bar T^3,\bar\varphi_1)$.
(To see this explicitly, note that $\varphi_1$ assigns the standard frame
$i,j,k$ in $\R^3$ to each point $(x,y,z) \in T^3 = \R^3/\Z^3$, and the
reflection $(x,y,z) \mapsto (x,y,-z)$ carries this frame onto $i,j,-k$.)
It follows by \ref{thm:d}d and \ref{thm:h}d that the total defect
vanishes.\footnote
{There is another way see this:  $\varphi_1$
obviously has degree zero, since it is an honest framing, and its
Hirzebruch defect can be computed by viewing $(T^3,\varphi_1)$ as the
boundary of the spin $4$-manifold $N$ obtained by removing an open tubular
neighborhood of a regular fiber in the rational elliptic surface of
signature $-8$ and Euler class $12$.  (This manifold is discussed in
detail in Chapter V of \cite{kir}.)  Indeed it can be shown that
$p_1(N,\varphi_1) = -2\chi(N)$, and so $h(\varphi_1) = -2\chi(N) -
3\sigma(N) = 0$.}
With each of the other seven spin structures, $T^3$ is
diffeomorphic to the boundary of the spin $4$-manifold $W = T^2\times B^2$
(with the Lie spin structure on $T^2$) -- note that the diffeomorphism
might not be the obvious one.   Now there is a natural framing on $W$,
namely the product of the Lie framing on $T^2$ with the constant framing
on $B^2$.  The restriction of this framing to the boundary is a stable
framing $\phi_0$ of $T^3$ with vanishing total defect, since $\chi(W) =
\sigma(W) = 0$.  In particular $\phi_0$ is the stabilization of the honest
framing $\varphi_0$ which assigns the frame $\cos(2\pi z)\,i + \cos(2\pi
z)\,j + k$ to the point $(x,y,z) \in T^3$.

Summarizing, we have shown that
$$
H(T^3,\varphi_0) = H(T^3,\varphi_1) = (0,0).
$$
This identifies the unique canonical framings in the eight spin structures
on $T^3$, since $\varphi_0$ represents each of the seven non-Lie
structures under a suitable diffeomorphism.  In particular, note that
$\lambda(\Sigma) = 0$ for every spin structure on $T^3$ (see Figure 2a),
which is consistent with Theorem \ref{thm:l} since $r(T^3) = 3$, and
$\mu(\Sigma) = 8$ for the Lie spin structures and $0$ for the rest.
\end{ex}

\begin{ex}\label{ex:products}
(Products) \ Let $M_g = F \times S^1$, where $F$ is a closed
orientable surface of genus $g$.  Then $M_g$ has $2^{2g+1} = |H^1(M;\Z_2)|$
spin structures.  Up to a diffeomorphism, however, there are only
four when $g>1$, two when $g=1$, and one when $g=0$.  To see this, assign
to each spin structure $\Sigma$ on $M_g$ the bordism invariant
$$
(\alpha,\beta) \in \Omega_2^{\text{spin}} \oplus
\Omega_1^{\text{spin}} = \Z_2 \oplus \Z_2
$$
determined by the restrictions of $\Sigma$ to the two factors of $F \times
S^1$ ($\alpha$ is the Arf invariant of $\Sigma|F$; see Chapter IV in
\cite{kir}).  Then an analysis of the action of the diffeomorphism group of
$M_g$ on $H^1(M_g;\Z_2)$ shows that any other spin structure with the same
bordism invariant is equivalent to $\Sigma$, and conversely for $g>1$, any
spin structure equivalent to $\Sigma$ will have the same bordism
invariant.  For $g = 0$ or $1$, the internal symmetry of the manifold $M_g$
further restricts the number of diffeomorphism classes of spin structures.
We consider only the cases $g=0$ and $g>1$, since $M_1 = T^3$ was
treated in the last example.  (See also Theorem \ref{thm:circlebundles}
for certain circle bundles over $F$.)

First observe by Theorem \ref{thm:l} that $\lambda(\Sigma) = 0$ for every
spin structure $\Sigma$ on $M_g$.  Indeed $r(M_g) = 2g+1$, and so this is
equivalent to the fact that $\mu(\Sigma) \equiv 0 \pmod4$.  In fact an
elementary induction shows that $\mu(\Sigma) \equiv 0 \pmod8$, since
$M_g$ can be viewed as a ``fiber connected sum" of $M_{g-1}$ and $M_1$
(i.e.\ remove tubular neighborhoods of circle fibers in both and identify
the boundaries) and the mu invariants clearly add.  It follows that the
total defect of any canonical stable framing on $M_g$ is $(0,0)$.

\smallskip
{\bf a)} For $g=0$ the manifold is $M_0 = S^2\times S^1$, which has two
spin structures with bordism invariants $(\alpha,\beta) = (0,0)$ and
$(0,1)$. In fact these spin structures are equivalent under the
automorphism $\tau$ of $S^2\times S^1$ which spins the $S^2$ once along
the $S^1$ factor.  Thus it suffices to consider the spin structure with
bordism invariant $(0,1)$, for which there is an obvious stable framing
$\phi_1$ obtained by restricting the product framing on $B^3 \times S^1$
(the constant framing on $B^3$ crossed with the tangent framing on $S^1$)
to the boundary.  Since the Euler characteristic and signature of $B^3
\times S^1$ vanish, the total defect of this framing vanishes as well.
Therefore $\phi_0 = \tau^*\phi$ and $\phi_1$ are the unique canonical
framings for the two spin structures on $S^2\times S^1$.

Note that since the stable framings $\phi_i$ have degree zero, they are
stabilizations of honest framings $\varphi_i$ on $S^2\times S^1$, obtained
from $\phi_i$ as follows:  Rotate each frame so that the last vector, which
is initially tangent to $S^1$, becomes the inward normal to $B^3$, and
then drop this vector from the frame.

\smallskip
{\bf b)} For $g>1$ we only treat the case when the Arf invariant $\alpha
= 0$.  Then $M_g$ is the spin boundary of a product $4$-manifold $N\times
S^1$.  Now if $\beta = 1$, then the canonical framing in the associated
spin structure $\Sigma_1$ is just the restriction $\phi_1$ to
$M_g$ of any product framing on $N\times S^1$.  The total defect of this
framing vanishes since the Euler characteristic and signature of $N\times
S^1$ are zero.  For the spin structure $\Sigma_0$ with $\beta = 0$, we
modify $\phi_1$ by putting a full twist in the plane of the first two
vectors in each frame while traversing the $S^1$ factor.  This yields the
canonical stable framing $\phi_0$ compatible with $\Sigma_0$.  Of course
both $\phi_0$ and $\phi_1$ can be homotoped to honest framings as in the
genus $0$ case.
\end{ex}

\heading{$2$-framings}

There is a canonical spin structure on $2\tau_M$ coming from the diagonal
embedding of $\tau_M$ in $2\tau_M$, namely $\Sigma\oplus\Sigma$ for any
spin structure $\Sigma$ on $M$.  It is easy to see that this spin
structure is in fact independent of the choice of $\Sigma$.  Following
Atiyah, we consider only the $2$-framings which extend this canonical spin
structure (although the general case is not much harder).   These
$2$-framings form an affine space $\tF$ with translation group $\pi_3(SO_6)
= \Z$.  In particular the generator $\sigma$ acts in the usual way on the
first four vectors in any $2$-framing $\tf$, and trivially on the last two,
and it follows that
$$
h(\tf+\sigma) = h(\tf)+2.
$$
Recalling from \S1 that $h(2\varphi)$ is even for any honest framing
$\varphi$ of $M$, we see that $\tF$ can be identified with the even
integers in $\Z = (h)$.  It follows that there exists a {\sl unique}
$2$-framing in $\tF$ with zero Hirzebruch defect.  This is Atiyah's {\it
canonical} $2$-framing.

\begin{defn}\label{defn:ctf}
(Atiyah) The {\it canonical} $2$-framing on $M$ is the unique $2$-framing
$\ctf$ compatible with the canonical spin structure on $2\tau_M$
for which $h(\ctf) = 0$.
\end{defn}

The canonical $2$-framing $\ctf$ on $M$ is clearly given by
$
\ctf = 2\varphi - h(\varphi)\sigma
$
for {\sl any} honest framing $\varphi$.  In particular $\ctf$ is of the
form $2\varphi$ for some framing on $M$ if and only if $\varphi$ is the
canonical framing in a spin structure $\Sigma$ with $\lambda(\Sigma) = 0$.
By Theorem \ref{thm:l}, this will occur if and only $M$ has a spin
structure $\Sigma$ with $\mu(\Sigma) \equiv 2(1+b_1(M))\pmod 4$, where
$b_1$ is the first Betti number.  This condition is satisfied, for
example, for all products $F\times S^1$ (see Example \ref{ex:products}),
but fails for all homology spheres.  It holds for the lens spaces $L(n,1)$
if and only if $n\equiv\sign(n)\pmod4$.

More generally, one can ask whether the canonical framing can be written
as a Whitney sum $\varphi_1\oplus\varphi_2$ of two framings (necessarily)
in the same spin structure $\Sigma$.  This is clearly equivalent to
having $\lambda(\Sigma) = 0$ or $2$, which by Theorem \ref{thm:l} occurs
if and only if the number $s(M)$ of $2$-primary components in $H_1(M)$ is
even.  For example this condition is satisfied for all homology spheres
(e.g.\ the canonical $2$-framing on $S^3$ is the sum of the two Hopf
framings) and all odd lens spaces, but fails for all even lens spaces.


\section{Other natural framings}

In this section we consider a variety of naturally arising framings of
$3$-manifolds, and show how to compute their total defects.  As in
previous examples, we may write $(M,\phi)$ for $(\phi)$ to highlight the
underlying manifold $M$ on which the framing $\phi$ is defined.

\heading{Homogeneous spaces}

Let $G$ be a finite subgroup of $S^3 = SU_2$.  Then the right handed Hopf
framing $\varphi_+$ on $S^3$, which is defined by left translation of a
frame at the identity, induces an honest framing on the homogenous space
$G\backslash S^3$ of right cosets of $G$ in $S^3$.  This framing
will also be denoted by $\varphi_+$, with associated spin structure
$\Sigma_+$, and its Hirzebruch defect can be computed from the formula in
Theorem \ref{thm:h}c to be
$$
h(G\backslash S^3,\varphi_+) = (2-\sigma(G))/|G|
$$
where $\sigma(G)$ is three times the signature defect of the universal
covering $S^3 \to G\backslash S^3$.  By Examples 3 and 4 in Appendix B,
$\sigma(G)$ is equal to $m^2-3m+2$, $4m^2+2$, $98$, $242$ or $722$,
according to whether
$G = C_m$ (for $m\ge1$), $D_{m}^*$ (for $m\ge2$), $T^*$, $O^*$ or $I^*$.
(The star indicates the double cover of the relevant subgroup of $SO_3$
\cite{wolf}.)  The corresponding Hirzebruch defects are $3-m$, $-m$, $-4$,
$-5$ and $-6$, respectively.  Note that $\varphi_+$ is canonical only in
the cyclic case for $m\le5$, and in the dihedral case for $m = 2$
when $G$ is the quaternion group $Q_8$.

In particular $C_m\backslash S^3$ is the lens space $L(m,1)$ (for $m>0$).
The quotient framing $\varphi_+$ with Hirzebruch defect
$$
h(L(m,1),\varphi_+) = 3-m
$$
can then be modified to give a canonical framing $\varphi_+ + \lfloor
(m-1)/4 \rfloor \rho$ for the associated spin structure $\Sigma_+$ on
$L(m,1)$.  If $m$ is even, then $L(m,1)$ has an another spin structure,
and it will be seen below how to construct an associated framing.\footnote
{The case $m=2$ has already been fully treated in \ref{ex:lie}, since
$L(2,1) \ (=\R P^3) = SO_3$.  Moreover, since $h(\varphi_+) = 1$, the
quotient framing $\varphi_+$ coincides with the right-handed Lie framing.
Similarly the quotient framing $\varphi_-$ for the left coset space
$S^3/C_2$ coincides with the left handed Lie framing.}

For the Poincar\'e homology sphere $P^3 =I^*\backslash S^3$ we have
$$
h(P^3,\varphi_+) = -6.
$$
Thus $\varphi_+ + \rho$ is a canonical framing (with defect $-2$) for the
unique spin structure on $P^3$.  The other canonical (stable) framings are
obtained by adding $\rho$, $\sigma$ or $\rho-\sigma$ to this one.

\heading{Circle bundles}

Consider an oriented circle bundle $\E$ with Euler class $n$ over a
closed, oriented surface $F$ of genus $g$.  Set $\chi = 2-2g$, the
Euler characteristic of $F$.  There is an obvious vector field $\tau$ on
$\E$, tangent to the oriented circle fibers.  Any framing of $E$ which
extends $\tau$ will be called {\it fiber-preserving}.  The following lemma
and discussion can be compared with Gompf's discussion of fields of
$2$-planes on 3-manifolds in section 4 of \cite{gompf}.

\begin{lem} $\E$ has fiber-preserving framings if and only if $n$ divides
$\chi$.
\end{lem}

\begin{proof}  Let $\tau^{\perp}$ denote the oriented $2$-plane bundle over
$\E$ which is orthogonal to $\tau$.  Clearly $\tau$ extends to a framing of
$E$ if and only if $\tau^\perp$ has a nonvanishing section.  The
obstruction to finding such a section is given by the Euler class
$e(\tau^\perp) \in H^2(E) = \Z^{2g} \oplus \Z_n$.  Since the projection
$p\!:\E \to F$ is covered by a bundle map $\tau^{\perp} \to \tau_F$, it
follows that $e(\tau^{\perp}) = p^*(e(\tau_F)) = (0,\chi)$, and so
$e(\tau^\perp) = 0$ if and only if $\chi \equiv 0 \pmod n$.
\end{proof}

The fiber-preserving framings of $E$ can be described explicitly as
follows.  Let $D$ be a disk in $F$ and let $F_0$ be the complement of the
interior of $D$.  Then $\E$ is trivial over $F_0$, so we frame its tangent
bundle by the product of $\tau$ with a (tangential) framing of $F_0$.  Note
that the framing of $F_0$ is not unique, for it can be changed by elements
of $H^1(F)$, but the $2$-frame on $\partial F_0$ is unique up to homotopy,
and it spins $1-\chi$ times compared to the stabilized tangent framing of
$\partial F_0$.  When the trivial circle bundle over $D$ is attached to the
trivial circle bundle over $F_0$ to get a bundle with Euler class $n$, we
pull back a framing on $\partial D \times S^1 = T^2$.  On an $(1,n)$ curve
in $T^2$ (the image of $\partial F_0 \times$point), the framing has one
vector equal to $\tau$ and the other two vectors spin $\chi-1$ times
compared to a tangent vector to the $(1,n)$ curve.  Thus, if we frame $D
\times S^1$ by $\tau$ and a pair of vectors which are constant on each
copy of $D$, but rotate $\chi/n$ times as the $S^1$ is traversed, then
this framing matches up with the one induced from the framing on $F_0
\times S^1$.

We illustrate this construction for the case $g=0$ and $n=1$, where the
fiber-preserving framing corresponds to the right handed Lie framing on
$S^3$, and show how to use it to calculate the first Pontrjagin number of
the disk bundle $\Delta$ associated with $E$.

\begin{ex} The boundary $S^3$ of the normal bundle $\Delta$ of $\C P^1$ in
$\C P^2$ is the Hopf circle bundle over $\C P^1$, of Euler class
$1$.  The right handed Hopf framing $\varphi_+$ on $S^3$ corresponds in
the above discussion to $\chi/n = 2/1 = 2$ rotations in the framing on $D
\times S^1$.  This gives a stable framing $\phi_+$ of $\tau_{S^3}$, the
outward normal plus $\tau$ for one complex line, and the framing of
$\tau^{\perp}$ for the other complex line.  This stable framing does not
extend over $\C P^1$, but it does provide a trivialization of the {\sl
complex} bundle $\tau_{\Delta}$ over the boundary $S^3$, and so
$$
p_1(\Delta,\phi_+) = c_1^2(\Delta,\phi_+) -2c_2(\Delta,\phi_+)
$$
(see Appendix A).  Now $c_1$ is the obstruction to extending this
trivialization over the 2-skeleton of $\Delta$, that is, over a fiber of
$\Delta$.  There is one full twist in the complex line spanned by the
outward normal and $\tau$, and two full twists in the orthogonal complex
line as noted above, so $3\,\C P^1$ is the Poincar\'e dual to $c_1$.  Thus
$c_1^2$ is the self intersection of $3\,\C P^1$.  Of course $c_2$
is the Euler class, so $p_1(\Delta, \phi_+) = 9 - 4 = 5$.

Note that this calculation is consistent with previous calculations of
the Hirzebruch defect: $h(\phi_+) = p_1(\Delta,\phi_+) - 3\sigma(\Delta) =
5-3 = 2$.  It also gives direct confirmation of the Hirzebruch signature
formula for $\C P^2$.  Indeed, for the remaining $4$-handle of
$\C P^2$ we calculate $p_1(\bar{B}^4, \bar{\phi}_+) = - 2$, so $p_1(\C
P^2) = 5 - 2 = 3$.
\end{ex}

Now in general, observe that the flexibility in the initial choice of the
framing of $F_0$ in the construction above, arising from the action of
$H^1(F)$, leads to a fiber-preserving framing in each spin structure on $E$
which is the pull back $p^*\Sigma$ of some spin structure $\Sigma$ on
$F$.  These will be called the {\it bundle spin structures} on $E$.

\begin{thm}\label{thm:circlebundles}
If $n$ divides $\chi$, then there is a unique fiber-preserving framing
$\varphi$ in each bundle spin structure on $E$, and $h(\varphi) = n +
\chi^2/n - 3\sign(n)$.
\end{thm}

\begin{proof}
The existence of $\varphi$ was shown above, and the uniqueness (up
to not necessarily fiber-preserving homotopy) is clear.  Calculating the
relative Pontrjagin number for the disk bundle $\Delta$ over $F$ bounded
by $E$, as in the example, we have
$$
p_1(\Delta,\varphi) = (1+\chi/n)F \cdot (1+\chi/n)F - 2\chi =
(1+\chi/n)^2n-2\chi = n + \chi^2/n.
$$
The Hirzebruch defect is then gotten by subracting $3\sigma(W) =
3\sign(n)$.
\end{proof}


\section{Boundaries and Surgery}

In this section we discuss natural stable framings of a $3$-manifold
$M$ that arise when viewing it as the boundary of a simply-connected spin
$4$-manifold $W$.  For simplicity, we consider only the case when $W$ can
be built without $1$- or $3$-handles, and so $M$ can also be viewed as
surgery on a link in $S^3$.

\heading{Boundaries}

Let $M_L$ denote the boundary of the $4$-manifold $W_L$ obtained by adding
$2$-handles to $B^4$ along a (normally) framed link $L$ in $S^3$
\cite{kir}.  It is a classical theorem of Lickorish and Wallace
that every closed connected oriented $3$-manifold is diffeomorphic to some
$M_L$.

If the framings on $L$ are all even, then $W_L$ is parallelizable with a
unique framing (since $H_1(W_L) = 0$). This framing restricts to a stable
framing $\d_L$ of $M_L$, with associated spin structure $\Sigma_L$.
By Theorems \ref{thm:d}b and \ref{thm:h}b, the total defect $H = (d,h)$ of
this stable framing is
$$
H(\d_L) = (\chi_L,-3\sigma_L)
$$
where $\chi_L$ and $\sigma_L$ denote the Euler characteristic and
signature of $W_L$.  From $\d_L$ we get an honest framing $\e_L =
\d_L + \chi_L\,\sigma$ with Hirzebruch defect
$$
h(\e_L) = \lambda_L = 2\chi_L-3\sigma_L.
$$
Note that $\lambda_L$ reduces mod $4$ to the invariant $\lambda(\Sigma_L)$
that identifies the affine lattice in the $dh$-plane associated with the
set $\F_L$ of stable framings compatible with $\Sigma_L$ (see above
Theorem \ref{thm:l}).  Now subtracting a suitable multiple of $\rho$ from
$\e_L$ gives a canonical framing for $\Sigma_L$ with total defect
$(0,\lambda)$, for $|\lambda|\le2$ and $\lambda \equiv \lambda_L$.  Since
every spin $3$-manifold $(M,\Sigma)$ is spin diffeomorphic to some
$(M_L,\Sigma_L)$, and there is an effective algorithm for finding $L$
\cite{kap}, this provides a general construction for all the canonical
framings of any given $3$-manifold.

\penalty-800
\begin{ex}
The lens space $L(m,1)$ for $m>0$ can be described either as the boundary
of $W_K$ where $K$ is the unknot with framing $-m$ (the minus sign
can be seen geometrically from a careful look at the quaternions), or as
the boundary of $W_L$ where $L$ consists of a simple chain of
$m-1$ simply linked unknots all with framings $+2$.  Now $L(m,1)$ has
exactly two spin structures if $m$ is even, given by $\Sigma_K$ and
$\Sigma_L$, and a unique spin structure if $m$ is odd, given by
$\Sigma_L$, with mu invariants $\mu_K \equiv -1$ and $\mu_L \equiv
m-1$.\footnote {This is a standard exercise in the calculus of framed
links:  Each spin structure $\Sigma$ on a $3$-manifold $M_L$ is uniquely
specified by a {\sl characteristic} sublink $C$ of $L$, with $\mu(\Sigma)
\equiv \sigma_L - C\cdot C + 8\Arf(C)\pmod{16}$ (see \cite[Appendix
C]{kme}).  The algorithm in \cite{kap} shows how to {\sl blow down} $C$ to
obtain an even framed link representing $\Sigma$ as above.}
The corresponding stable framings $\d_K$ (for even $m$) and $\d_L$ (in
general) have total defects
$$
H(L(m,1),\d_K) = (2,3) \qquad\text{and}\qquad H(L(m,1),\d_L) =
(m,3-3m)
$$
since $\chi_K = 2$, $\sigma_K = -1$, $\chi_L = m$ and $\sigma_L = m-1$.
From this (or using Theorem \ref{thm:l} and the mu invariant calculation
above) we deduce that $\lambda_K \equiv -1$ and $\lambda_L \equiv 3-m
\pmod4$.  The associated canonical framings are $\d_K+2\sigma-2\rho$ and
$\d_L+m\sigma+{1\over4}(m-3+\lambda)\rho$, for $|\lambda|<2$ and
$\lambda\equiv 3-m\pmod4$.

The stable framing $\d_L$ can also be expressed in terms of the quotient
$\varphi_+$ of the right-handed Lie framing on $S^3$, discussed in the
previous section.  Indeed the defect calculations above show that $\d_L =
\phi_+ - m\sigma$, where $\phi_+$ is the stabilization
$\nu\oplus\varphi_+$.  Implicit in this statement is the assumption that
$\phi_+$ is compatible with the spin structure $\Sigma_L$.  This is
automatic if $m$ is odd, when there is only one spin structure on
$L(m,1)$, and is also clear for $m
\equiv 2\pmod4$, since $\lambda(\Sigma_L)$ then distinguishes the two spin
structures.

If $m\equiv0\pmod4$, then the situation is more subtle.  In this case the
affine lattices for $\F_K$ and $\F_L$ coincide, both being equal to
$\Lambda_{-1}$ (see \S2).  It follows from Theorems \ref{thm:d} and
\ref{thm:h} applied to the universal cover $\pi\!:S^3\to L(m,1)$ that
$\pi^*$ embeds both lattices onto the same sublattice $m\Lambda_0 +
(0,2)$ of $\Lambda_2$ (of index $m^2$).  This implies that there exist
{\sl homotopic} representatives of the Lie framing $\varphi_+$ on
$S^3$, equivariant with respect to the covering transformations of $\pi$,
which project to framings $\varphi_+^K$ and $\varphi_+^L$ of defect $3-m$
in the respective spin structures $\Sigma_K$ and $\Sigma_L$ on $L(m,1)$.
Thus the precise statement when $m$ is a multiple of $4$ is that $\d_K =
\phi_+^K + (1+m/4)\rho - 2\sigma$ and $\d_L = \phi_+^L - m\sigma$.
\end{ex}

\begin{ex}
The Poincar\'e homology sphere $P^3$ is the boundary
of $W_L$ where $L$ is the $E_8$-link with framings $+2$.  Since
$\chi_L = 9$ and $\sigma_L = 8$, we have
$$
H(P^3,\d_L) = (9,-24)
$$
and so $\d_L + 9\sigma + r\rho$ is canonical for $r=1$ or $2$.  Also note
that $\d_L + 9\sigma = \phi_+$, the stabilized quotient of the right-handed
Lie framing on $S^3$.
\end{ex}

\penalty-800
\heading{Surgery}

Starting with a framed link $L$ in $S^3$ as above, there is a natural way
to construct a $2$-framing $\tf_L$ of $M_L = \partial W_L$, discovered by
Freed and Gompf \cite{fg}, which proceeds by extending the canonical
$2$-framing on the complement of $L$ over the solid tori which are glued
in under surgery on $L$. They showed that
$$
h(\tf_L) = 2\tau_L -6\sigma_L
$$
where $\tau_L = \sum a_i$, the sum of the framings $a_i$ of the components
of $L$, and $\sigma_L = \sigma(W_L)$.

Using the same philosophy, we show how to construct framings and
stable framings whose Hirzebruch defects satisfy similar formulas.  For
simplicity we assume that the (normal) framings on the components of $L$
are all {\sl even}, and consider only framings of $M_L$ which are
compatible with the spin structure $\Sigma_L$coming from $W_L$.  As above,
$\chi_L$ will denote the Euler characteristic of $W_L$.  Note that $\chi_L
= \ell+1$, where $\ell$ is the number of components in $L$.

Our constructions are based on the following gluing principle for
$4$-manifolds.

\begin{lem} {\rm (Gluing)}\label{lem:glue}  Let $W_1$ and $W_2$ be compact
oriented $4$-manifolds with stable framings $\phi_1$ and $\phi_2$ on the
boundary, and $f\!:M_2\to M_1$ be an orientation reversing diffeomorphism
between codimension zero submanifolds $M_i \subset \partial W_i$ such that
$f^*\phi_1 = \bar\phi_2$.  Then there is a natural stable framing $\phi$ on
the boundary of the $4$-manifold $W = W_1\cup_f W_2$ such that  $d(\phi) =
\sum d(\phi_i) - \chi$, where $\chi$ is the Euler characteristic of
(either) $M_i$, and $p_1(W,\phi) = \sum p_1(W_i,\phi_i)$.\footnote
{The behavior of the signature under gluing, discussed in
\cite{wall}, reduces to the computation of $\sigma(\Theta)$ since
$\sigma(W) = \sum\sigma(W_i) + \sigma(\Theta)$ by Novikov additivity.}
\end{lem}

\begin{proof}
After a homotopy, we can assume that $\phi_1|M_1$ and $\bar\phi_2|\bar M_2$
are identified under $f$.  Extend each $\phi_i$ to a framing of a collar
neighborhood of $\partial W_i$ inside $W_i$.  The result of gluing $W_1$
and $W_2$ together using $f$ is a priori a $4$-manifold with {\sl corners},
which we denote by $V$, and $\phi_1$ and $\phi_2$ clearly combine to give
a framing $\Phi$ of a collar neighborhhood of $\partial V$.
Now $W$ can be viewed as a retracted copy of $V$ with smooth boundary in
this collar, and the restriction of $\Phi$ to $\partial W$ is
the desired stable framing $\phi$.

Note that $W$ is the union of (copies of) $W_1$, $W_2$, and a ``thin"
$4$-manifold $\Theta$ bounded by $\partial W \cup \partial \bar W_1
\cup \partial \bar W_2$.  By Theorem \ref{thm:d}b $d(\phi) - d(\phi_1) -
d(\phi_2) = \chi(\Theta) = -\chi$, and the additivity of $p_1$ is obvious.
\end{proof}

Now apply the lemma to the handlebody $W_L$, obtained from $B^4$ (the
$0$-handle) by attaching $2$-handles along the even framed link $L$.
Note that the attaching regions are solid tori, with $\chi = 0$, and so
the degrees add.  We begin with the stable case.

\begin{con} (Natural stable framings of $M_L$ extending $\Sigma_L$.)
Because the framings on $L$ are even, any choice of stable framing on
$S^3 = \partial B^4$ will extend across the $2$-handles.  If we choose
$\d+n\sigma$ (where $\d$ is the canonical stable framing coming from
$B^4$) and frame each $2$-handle with $\d$, then the associated stable
framing of $M_L$, denoted $\phi_L^n$, will have total defect
$$
H(\phi_L^n) = (\chi_L-n,2n-3\sigma_L).
$$
In particular the case $n=0$ gives the restriction $\d_L$ of the
unique framing on $W_L$, and $n=\chi_L$ gives the honest framing
$\e_L$ discussed above.

Another natural choice is $n=\tau_L/2 = \sum b_i$ (where the even framings
on the components of $L$ are $a_i = 2b_i$), and we write $\phi_L$ for the
corresponding stable framing.  This stable framing, with Hirzebruch defect
$$
h(\phi_L) = \tau_L - 3\sigma_L
$$
reminiscent of the Freed-Gompf $2$-framing, arises in an effort to constuct
an ``explicit" framing of $M_L$, one that can be visualized without
appealing to homotopy theory.  The construction is in two steps:

\medskip
\begin{enumerate}
\item[$\bullet$] {\sl Tilting:} Choose a framing representing $\d$ which
reflects the topology of $L$.\smallskip
\item[$\bullet$] {\sl Matching:} Adjust the framing to take the
framings on $L$ into account.
\end{enumerate}
\medskip

{\sl Tilting.} \ For simplicity we assume that $L$ consists of
a single knot $K$; the process that we describe can be
repeated for each component of $L$.  Start with an oriented embedding of
$K$ in a coordinate chart $\R^3$ in which $\d$ is the standard constant
framing (outward normal ``$1$" followed by $i,j,k$), and assume that the
projection of $K$ onto the $ij$-plane is generic.  Let $w$ be the writhe
(i.e.\ the sum of the crossing signs) and $d$ be the Whitney degree (i.e.\
the degree of the Gauss map) of this projection.   Then there is a
canonical shortest isotopy of the framing near $K$ (called ``tilting")
which rotates $k$ to the oriented tangent vector of $K$.  The vectors $i$
and $j$ will then provide a (normal) framing of $K$ which spins $w-d$ times
relative to the $0$-framing.

Note that $w-d$ is unchanged by the second and third Reidemeister moves,
but can change by $\pm 2$ by the first move: when a kink (there are four
types depending on the sign of the crossing and the orientation) is
added, $w$ and $d$ both change by $\pm 1$, and so $w-d$ changes by $-2, 0$
or $2$.  Thus a projection of $K$ can be chosen so that the quantity $w-d$
takes on any prescribed {\sl odd} value (since it is odd for the round
unknot, and crossing changes preserves its parity).  We choose a
projection with $w-d = 1$.  This yields a stable framing on $S^3$ whose
second vector is tangent to $K$ and whose last two provide the
$+1$-framing of $K$ in $S^3$

\medskip

{\sl Matching.} \ The framing on the $2$-handle $H = B^2\times
B^2$ can be taken to be the product of the constant framings on the
two factors, denoted by $u_1, v_1$ and $u_2, v_2$.  It is more
useful however to have the framing on the attaching circle $K = S^1\times
B^2$ consist of one vector tangent to $K$ and one pointing into $H$. This
can be achieved by an isotopy of the framing on $H$ with support near $K$
which rotates $u_1$ to the oriented tangent vector to $K$, rotates $v_1$ to
point inward, and necessarily moves $u_2$ and $v_2$ so that they rotate
$-1$ time as an attaching circle is traversed. (The last statement follows
from the well known fact that the natural homomorphism $\pi_1(U_1) \oplus
\pi_1(U_1) \rightarrow \pi_1(U_2)$ maps $(p,q)$ to $p+q$.)

First suppose that $K$ has the $0$-framing.  Then the framing on $K$ in
the $0$-handle perfectly matches its framing in the $2$-handle, with
$\{u_1, u_2, v_1, v_2\}$ corresponding in order to $\{i, 1, j,k\}$; note
the transposition of $1$ and $i$ which corresponds to the fact that one
orientation is negated when two manifolds are joined along part of their
boundaries.  Also note that the framing along $K$ is $+1$ compared to the
$0$-framing of $K$, whereas the normal bundle of the attaching circle
has framing $-1$ in the boundary of the $2$-handle; this is correct because
of the orientation change.\footnote
{The easiest way to be convinced of this is to consider a Hopf link in
$S^3$, which is the equator of $S^4$.  Orient the Hopf link, and let it
bound $2$-balls in both hemispheres of $S^4$. The two $2$-balls intersect
$+1$ in one hemisphere and $-1$ in the other, since the two $2$-spheres
intersect algebraically zero, and so the linking number of the oriented
Hopf link is $+1$ in one hemisphere and $-1$ in the other.}
Since the framings agree on $K$, they are close on a neighborhood of $K$
and therefore canonically isotopic on this neighborhood, and so the gluing
lemma applies.

Now suppose that $K$ has framing $a = 2b$.  Recall from \S2 that adding
$\sigma$ to a stable framing puts a full right twist in both the normal
and conormal planes along each diameter of a small $3$-ball.  In fact this
local description can be modified by a homotopy so that along one diameter,
the normal plane rotates two full right twists while the conormal plane
does not rotate at all.  Using the latter description,  modify the stable
framing $\d$ near $K$ by adding $b = (\tau_L/2)$ copies of $\sigma$ in a
ball intersecting $K$ in a diameter.  This puts $2b = a$ twists in the
framing along $K$ in the $0$-handle, and so the framings along $K$ in the
$0$- and $2$-handles match.
\end{con}

\begin{con} (Natural framings of $M_L$ extending $\Sigma_L$.)
We apply Lemma \ref{lem:glue} as in the stable case, framing the boundary
of the $0$-handle with  $\varphi_\pm + n\rho$ (where $\varphi_\pm$ are the
canonical Lie framings) and each $2$-handle with $\varphi_\pm$.
The associated framings $\varphi_{\pm L}^n$ have Hirzebruch defects
$$
h(\varphi_{\pm L}^n) = 4n \pm 2\chi_L - 3\sigma_L
$$
and so in particular $\varphi_{+L}^0 = \e_L$.  As for stable framings, the
case $n=\tau_L/2 = \sum b_i$ (where the $a_i = 2b_i$ are the framings on
$L$) is also a natural choice, and we write $\varphi_{\pm L}$ for the
corresponding framings of $M_L$.  These framings have Hirzebruch defects
$$
h(\varphi_{\pm L}) = 2\tau_L - 3\sigma_L \pm 2\chi_L
$$
and can be constructed by a ``tilting/matching" scheme as above.  The
only difference is that when it becomes necessary to modify the framing
near $L$ in $S^3$ to match the $2$-handle framings, we use $\rho$ (which
also puts two full twists in the framing along a diameter) instead of
$\sigma$.  Note that $\tau_L$ is replaced with $2\tau_L$, since adding
$\rho$ adds $4$ to $p_1$, and thus to $h$, rather than $2$.
\end{con}

\begin{con} (Natural $2$-framings on $M_L$.)
Finally for $2$-framings, one way to proceed is to take the Whitney sum of
a pair of natural framings (just described) on each factor of $\tau_M$.
In particular, the natural $2$-framing of Freed and Gompf \cite{fg} can be
expressed as
$$
\tf_L = \varphi_{+L} \oplus \varphi_{-L}^0 =
\varphi_{+L}^0 \oplus \varphi_{-L}.
$$

Alternatively $\tf_L$ can be described in a natural way using an
analogue of Lemma \ref{lem:glue}.   For the reader's convenience, we also
recall the description of $\tf_L$ in \cite{fg}:  The first step is to
isotope the canonical $2$-framing on $S^3$ so that restricted to the
boundary torus $T$ of the attaching region $V$ for each $2$-handle, it is
the Whitney sum of two copies of a Lie framing on $T$ plus the normal
vector.  There is a choice in how this is done, since $\pi_3(\Spin_6) =
\Z$, but a different choice would change the $2$-framing on $S^3 - V$ by
some $\alpha \in \pi_3(\Spin_6)$ while changing the $2$-framing on $V$ by
$-\alpha$, and so these changes cancel on $M_L$ after regluing.  Now the
solid torus is removed and glued back in by an element $A$ of $SL_2(\Z)$,
but the $2$-framings on $T$ can be made to match by choosing a shortest
path in $SL_2(\R)$ from $A$ to the identity.  This gives $\tf_L$, and Freed
and Gompf prove directly that $h(\tf_L) = 2\tau_L - 6\sigma_L$.
\end{con}


\appendix


\section{Pontrjagin numbers}

This appendix is a brief review of some aspects of the theory of
characteristic classes from the obstruction point of view.  The focus is on
the relative first Pontrjagin number of a compact $4$-manifold, which
along with the signature is used to define the Hirzebruch defect of a
framing of the boundary $3$-manifold.  This material is of course well
known, although we do not know where to find an elementary discussion of
the relative theory from the obstruction point of view.

\heading{Absolute characteristic classes}

Let $X$ be a closed oriented $4$-manifold.  The $k$th Chern class
$c_k(\omega)$ of a complex $n$-plane bundle $\omega$ over $X$ can be
identified with the obstruction to finding an {\it $(n-k+1)$-field} on
$\omega$ (i.e.\ $n-k+1$ linearly independent sections of $\omega$) over
the $2k$-skeleton of $X$ \cite[p.171]{ms}.  In particular the second Chern
class $c_2(\omega) \in H^4(X;\Z)$ is the obstruction to finding an $(n\! -
\!\! 1)$-field on $\omega$ over all of $X$.  (This obstruction is a priori
an element of $H^4(X; \{ \pi_3(U_n/U_1) \})$, but the coefficient groups
are canonically identified with $\pi_3(SU_n) = \Z$ as explained in \S1.)
Evaluating on the fundamental class $[X]$ gives the second Chern number of
$\omega$, an integer which by abuse of notation will also be denoted by
$c_2(\omega)$.  For example, if $n=2$ then $c_2(\omega) = e(\omega_\R)$,
the Euler class of the underlying oriented real bundle or equivalently the
self intersection of the zero-section of $\omega$.  In particular
$c_2(\tau_X) = \chi(X)$ if $X$ is a complex surface, where $\chi$ is the
Euler characteristic.

The first Pontrjagin class of a real $n$-plane bundle $\xi$ over $X$
is defined in terms of the complexified bundle $\xi_\C$ by
$$
p_1(\xi) = -c_2(\xi_\C)
$$
\cite[p.174]{ms}.  If $\xi$ admits a complex structure, then $p_1$ can also
be computed using the formula $p_1(\omega_\R) = c_1^2(\omega) -
2c_2(\omega)$ \cite[p.177]{ms}.  For example $p_1(\tau_X) = -2\chi(X)$ if
$X$ is a complex elliptic surface, since $c_1^2=0$ for elliptic surfaces
\cite{kod}.

Taking $\xi=\tau_X$ and evaluating on the fundamental class gives
the first Pontrjagin number
$$
p_1(X) = p_1(\tau_X)[X].
$$
Thus $p_1(X)$ is seen to be the integer obstruction to the existence of a
$3$-field on the complexified tangent bundle of $X$.

\heading{Relative characteristic classes}

Let $W$ be a compact oriented $4$-manifold with nonempty boundary.  Then of
course $c_2(\omega)$ and $p_1(\xi)$ vanish for all complex bundles $\omega$
and real bundles $\xi$ over $W$, since $H^4(W) = 0$, and so we consider
{\sl relative} characteristic classes instead.  In particular for any
$(n\! - \!\! 1)$-field $\phi$ on $\omega|\partial W$, the relative second
Chern class $c_2(\omega,\phi) \in H^4(W,\partial W;\Z)$ can be defined as
the obstruction to extending $\phi$ over $W$.  One then defines the
relative first Pontrjagin class
$$
p_1(\xi,\phi) = -c_2(\xi_\C,\phi_\C),
$$
where $\phi$ is an $(n\! - \!\! 1)$-field on $\xi|\partial W$ and
$\phi_\C$ is the induced $(n\! - \!\! 1)$-field on $\xi_\C|\partial W$.
For complex bundles $\omega$ one also has $p_1(\omega_\R,\phi) =
c_1^2(\omega,\phi) - 2c_2(\omega,\phi)$, as in the absolute case.

Taking $\xi=\tau_W$ and evaluating on the fundamental class gives a
relative first Pontrjagin number
$$
p_1(W,\phi) = p_1(\tau_W,\phi) [W,\partial W]
$$
associated to any given tangential $3$-field $\phi$ over $\partial W$.
Thus $p_1(W,\phi)$ is an integer invariant which measures the obstruction
to extending $\phi$ to a $3$-field on the complexified tangent bundle of
$W$.  If $\phi$ is a stable framing, then the notation $p_1(W,\phi)$ should
be interpreted to mean $p_1(W,\varphi)$, where $\varphi$ is obtained by
dropping the first vector of $\phi$.  (Note that $\varphi$ determines
$\phi$, using the orientation of $W$.)

Similarly, taking $\xi = 2\tau_W$, one defines $p_1(W,\phi) =
p_1(2\tau_W,\phi)[W,\partial W]$ for any $7$-field $\phi$ on
$\tau_W \oplus \tau_W$ over $\partial W$.  For example, $\phi$ can be taken
to be of the form $\nu\oplus\tf$, where $\nu$ is the outward normal vector
field in the first factor of $\tau_W$ and $\tf$ is a $2$-framing of
$\partial W$.  In this case, this invariant is also denoted $p_1(W,\tf)$.
Note that
$
p_1(W,\varphi_1\oplus\varphi_2) = p_1(W,\varphi_1) +
p_1(W,\varphi_2)
$
for any two framings $\varphi_1$ and $\varphi_2$ of
$\partial W$, by the product formula for Pontrjagin classes
\cite[p.175]{ms}.


\section{Signature theorems}

This appendix contains a discussion of the signature theorem for closed
$4$-manifolds (classically due to Hirzebruch, and in its equivariant
form, to Atiyah and Singer) and various ``defect" invariants of
$3$-manifolds that result from the failure of this theorem for
$4$-manifolds with boundary.  In particular, we give a formula which
relates the behavior of the Hirzebruch defect of framed $3$-manifolds under
a covering projection, with the signature defect of the covering (see Lemma
B.1).

\heading{The Hirzebruch signature theorem}

The first Pontrjagin number of a closed oriented $4$-manifold $X$ is
related to the signature $\sigma(X)$ of $X$ (the difference of the number
of positive and negative eigenvalues of the intersection form on $H_2(X)$)
by  Hirzebruch's signature formula
$$
p_1(X) = 3\sigma(X)
$$
\cite[p.86]{hir} (see also \cite{hirze, ms, kir}).  This obviously
fails for manifolds with boundary, since then $p_1$ is always zero while
$\sigma$ need not be.  In fact it even fails when $p_1$ is replaced with
the appropriate relative version, as will be seen below (cf.\ \cite{aps}).
Note that replacing $\tau_X$ by $2\tau_X$ gives $p_1(2\tau_X)[X] = 2p_1(X)
= 6\sigma(X)$ (see \cite[p.175]{ms}).

One important consequence of the signature theorem is the multiplicativity
of the signature under finite covers: if $\tilde X\to X$ is an $r$-fold
covering map of closed oriented $4$-manifolds, then
$$
\sigma(\tilde X) = r\,\sigma(X).
$$
This follows from the multiplicativity of $p_1$, which is evident from the
obstruction point of view.  Again this fails for bounded manifolds.

\heading{Defects}

The failure of the signature theorem for bounded $4$-manifolds gives rise
to an integer invariant of framed $3$-manifolds (cf.\ \cite{atiyah,aps}).
This is what we have called the Hirzebruch defect in \S1, and is
the key invariant used in defining canonical framings in \S2.  We recall
the definition from \S1:  The {\it Hirzebruch defect} of a framing or
stable framing $\phi$ of a closed oriented $3$-manifold $M$ is defined by
$$
h(\phi) = p_1(W,\phi) - 3\sigma(W)
$$
for any compact oriented $4$-manifold $W$ with oriented boundary $M$.
(This is seen to be independent of the choice of $W$ by Novikov
additivity of the signature and the signature formula for closed
manifolds.)  Similarly define
$$
h(\tf) = p_1(W,\tf) - 6\sigma(W)
$$
for any $2$-framing $\tf$ of $M$, in light of the observation above that
$p_1(2\tau_X)[X] = 6\sigma(X)$ for closed manifolds.  Note that
$
h(\varphi_1\oplus\varphi_2) = h(\varphi_1) + h(\varphi_2)
$
for any two framings $\varphi_1$ and $\varphi_2$ of $M$, since the first
Pontrjagin numbers add (see Appendix A).

Similarly the failure of the multiplicativity of the signature for finite
covering spaces of bounded $4$-manifolds leads to an invariant for
$3$-dimensional coverings.  First define the {\it signature defect}
of an $r$-fold covering map $\Pi\!:\tilde W \to W$ of compact
oriented $4$-manifolds to be
$
\sigma(\Pi) = r\sigma(W) - \sigma(\tilde W).
$
Of course this vanishes if the manifolds are closed, but is in general
nonzero.  From this one defines the {\it signature defect} $\sigma(\pi)$ of
an $r$-fold covering map $\pi\!:\tilde M \to M$ of closed oriented
$3$-manifolds, as follows.  If $\pi$ is a regular cover, then it
corresponds to a finite index normal subgroup $H$ of $\pi_1(M)$, or
equivalently to a homomorphism $\pi_1(M) \to G = \pi_1(M)/H$.  Since
$3$-dimensional bordism over a finite group $G$ is torsion \cite{cf}, some
finite number $m$ of copies of $\pi$ bounds over $G$.  In other words,
there exists a finite covering $\Pi\!:\tilde W\to W$ of compact oriented
$4$-manifolds with $m\pi = \partial \Pi$ ($=$ the restriction
$\Pi|\partial\tilde W \to \partial W$).  Now set
$$
\sigma(\pi) \ = \ {1\over m}\sigma(\Pi) \ = \ {1\over
m}(r\sigma(W)-\sigma(\tilde W)).
$$
This is clearly well-defined by Novikov additivity and the multiplicativity
of the signature for closed manifolds.  If $\pi$ is irregular, then the
core $\hat H$ of $H$ (i.e.\ the intersection of all the conjugates of $H$
in $\pi_1(M)$) is normal and of finite index in both $\pi_1(M)$ and $H$,
and so corresponds to a pair of finite regular covers $\hat\pi\!:\hat M\to
M$ and $\tilde\pi\!:\hat M\to \tilde M$.  In this case the signature defect
is defined by $s\sigma(\pi) = \sigma(\hat\pi)-\sigma(\tilde\pi)$, where
$s$ is the index of $\hat H$ in $H$.

In the same spirit one can define the {\it Pontrjagin defect} $p_1(\Pi)$ of
an $r$-fold covering map $\Pi\!:(\tilde W,\tilde\phi) \to (W,\phi)$ of
compact bounded $4$-manifolds by
$
p_1(\Pi)  = rp_1(W,\phi) - p_1(\tilde W,\tilde\phi).
$
Here $\phi$ and $\tilde\phi$ are framings or stable framings on the
boundary which are assumed to be compatible (i.e.\  $\tilde \phi =
\Pi^*\phi$, the pull back of $\phi$ under the covering map).  It is
clear from the characterization of $p_1$ as an obstruction that this
invariant is in fact {\sl identically zero}.

Finally the {\it Hirzebruch defect} $h(\pi)$ of an
$r$-fold cover $\pi\!:(\tilde M,\tilde\phi) \to (M,\phi)$ of closed
$3$-manifolds, with compatible (stable) framings, is defined by
$$
h(\pi)  = rh(\phi) - h(\tilde M,\tilde\phi).
$$
Observe that for $\pi = \partial\Pi$, we have by definition $h(\pi) =
p_1(\Pi) - 3\sigma(\Pi)  = -3\sigma(\pi)$.   In fact the Hirzebruch and
signature defects are always related in this way.

\begin{lem}
Let $\pi\!:(\tilde M,\tilde\phi) \to (M,\phi)$ be an $r$-fold covering map
of closed oriented $3$-manifolds with compatible (stable) framings.  Then
$h(\pi) = -3\sigma(\pi)$.  In other words
$$
h(\tilde M,\tilde \phi) = r\,h(\phi) + 3\,\sigma(\pi)
$$
where $\sigma(\pi)$ is the signature defect of $\pi$. \qed
\end{lem}

\begin{proof}
If $\pi$ is regular, then some multiple of $\pi$ bounds, $m\pi =
\partial\Pi$.  It follows that $mh(\pi) = h(m\pi) = h(\partial\Pi) =
-3\sigma(\Pi)$.  Dividing by $m$ gives $h(\pi) = -3\sigma(\Pi)/m
= -3\sigma(\pi)$.  For irregular $\pi$ there exist regular covers $\hat\pi$
and $\tilde\pi$ of degrees $rs$ and $s$ such that $\hat\pi = \pi\tilde\pi$,
as discussed above.  Then by definition $sh(\pi) = h(\hat\pi)-h(\tilde\pi)
= 3(\sigma(\tilde\pi)-\sigma(\hat\pi))$ (since $\hat\pi$ and $\tilde\pi$
are regular).  Dividing by $s$ gives $h(\pi) = 3(\sigma(\tilde\pi) -
\sigma(\hat\pi))/s = -3\sigma(\pi)$.
\end{proof}

To exploit this result, one must find ways to compute the signature
defect.  The most powerful tool for this purpose is the equivariant
version of the signature theorem, due to Atiyah and Singer
\cite{as} (see also \cite{hirze} and \cite{gordon}).

\heading{The Atiyah-Singer $G$-signature theorem}

Let $G$ be a finite group of order $r$ acting effectively on a closed
oriented $4$-manifold $\tilde X$ by orientation preserving
diffeomorphisms.  Since $G$ is finite, the orbit space $X = \tilde X/G$ is
an oriented rational homology manifold with signature $\sigma(X)$ defined
in the usual way.  Furthermore, each element $g\in G$ has an associated
{\it $g$-signature} $\sigma(g)$ ($=\sign(g,\tilde X)$ in the
literature) given by the difference of the traces of $g$ acting on the
positive and negative definite subspaces of $H^2(\tilde X;\R)$ with the
cup product form.  An elementary representation theory argument shows that
the sum of all the $g$-signatures is equal to $r\sigma(X)$, and so
$$
r\sigma(X) - \sigma(\tilde X) = \sum_{g\ne1}
\sigma(g)
$$
since clearly $\sigma(1) = \sigma(\tilde X)$.

As one might suspect, this equality fails if $\tilde X$ is allowed
to have boundary.  Thus one obtains in the usual way an invariant for free
actions of $G$ on closed oriented $3$-manifolds $\tilde M$, or
equivalently covering spaces $\pi\!: \tilde M \to M = \tilde M/G$.
Namely, choose any compact oriented $4$-manifold $\tilde W$ bounded by
$\tilde M$ over which the action extends, and write $W$ for the orbit
space of this extended action.  Then define
$$
\sigma(\pi) = r\sigma(W) - \sigma(\tilde W) - \sum_{g\ne1}
\sigma(g).
$$
This invariant clearly coincides with the signature defect of $\pi$
defined above.

Now the $G$-signature theorem \cite{as} gives local formulas for the
$g$-signatures $\sigma(g)$ when $g\ne1$ in terms of the infinitesimal
action of $g$ on its fixed point set.  In particular, this fixed point set
consists of a finite union of points $x_i$ and connected surfaces $F_j$,
and if $g$ acts on the tangent space at $x_i$ by rotating a pair of
orthogonal planes through angles $\alpha_i$ and $\beta_i$, and on the
tangent space at each point on $F_j$ by rotating the normal plane to $F_j$
through an angle $\gamma_j$, then
$$
\sigma(g) = -\sum \cot(\alpha_i/2)\cot(\beta_i/2) + \sum F_j\cdot
F_j\,\csc^2(\gamma_j/2).
$$
These formulas provide a way to compute signature defects in many
situations.

For example Hirzebruch \cite{hirze} used this approach to compute
the defects of (the universal covers of) lens spaces in terms of classical
Dedekind sums (see also \cite{km}).   This was accomplished by viewing
$L(m,n)$ as the orbit space of a linear action of the cyclic group $C_m$
on $S^3 = \partial B^4$.  When $n=\pm1$ the picture is particularly simple
since $C_m$ can be taken to be a subgroup of $S^3$ acting by
multiplication.  We explain this case below, along with the slightly more
complicated case of the binary icosahedral group $I^*$ of order 120.  The
quotient $I^*\backslash S^3$ is the well known Poincar\'e homology sphere
$P^3$.

\begin{ex} The lens space $L(m,1)$ for $m>0$ has signature defect
$(m-1)(m-2)/3$.   To see this, view $S^3$ as the unit quaternions and
$L(m,1)$ as the homogeneous space $C_m\backslash S^3$ of right cosets of
any of the cyclic subgroups $C_m$ of $S^3$ of order $m$ (they are all
conjugate).  The action of $C_m$ by left multiplication extends by
coning to $W = B^4$, with quotient $\tilde W=$ the cone on $L(m,1)$.  The
signatures of $W$ and $\tilde W$ vanish, as both spaces are contractible,
and so the defect consists only of the contributions from the fixed
points.

Now each $u\ne1$ in $C_m$ has a unique fixed point at the origin, where the
action is given by left multiplication by $u$.  This rotates the plane
spanned by $1$ and $u$ through the angle $\arg(u)$ ($= 2k\pi/m$ for
some $k$) and the orthogonal plane through the {\sl same} angle (it would
be the opposite angle if the action were on the right), and so contributes
$\cot^2(\arg(u)/2)$ to the defect.  Since the elements of $C_m$ are equally
spaced on a great circle in $S^3$ through $\pm1$ the total contribution is
$$
\sum_{k=1}^{m-1} \cot^2(k\pi/m) = (m-1)(m-2)/3.
$$
(The closed form for the sum is classical, cf.\ \cite[p.19]{hirze}.)
\end{ex}

\begin{ex} The Poincar\'e homology sphere $P^3 = I^*\backslash S^3$ has
signature defect $722/3$.  Indeed, proceeding as in the previous example,
each $u\ne1$ in $C_p$ contributes $\cot^2(\arg(u)/2)$ to the defect.  Note
in particular that $-1$ makes no contribution, since $\cot(\pi/2) = 0$.
Now observe that the icosahedral group $I$ can be expressed as the union of
fifteen cyclic subgroups of order $2$ (given by $\pi$-rotations around the
fifteen axes through opposite edges of the icosahedron), ten of order $3$
(given by $2k\pi/3$-rotations around the ten axes through opposite faces),
and six of order $5$ (given by $2k\pi/5$-rotations around the six axes
between opposite vertices).  These lift to cover $I^*$ with fifteen cyclic
subgroups of order $4$, ten of order $6$, and six of order $10$.  Since
any two of these subgroups intersect in $\pm1$, the defect is just the sum
of their contributions.  From the previous example we know that a subgroup
of order $m$ contributes $(m-1)(m-2)/3$, and so the defect is $15\cdot6/3
+ 10\cdot20/3 + 6\cdot72/3 = 722/3$.

A similar calculation gives $98/3$ and $242/3$ for the defects of quotients
$T^*\backslash S^3$ and $O^*\backslash S^3$ by the binary tetrahedral and
octahedral groups, and $(4m^2+2)/3$ for the defects of the prism
manifolds $D_{m}^*\backslash S^3$.  (See for example \cite[\S2.6]{wolf}
for a discussion of the finite subgroups of $S^3$.)
\end{ex}

\end{document}